\nonstopmode
\documentclass[11pt]{article}

\usepackage{amsthm}
\usepackage{amssymb}
\usepackage{amsmath}
\usepackage{bbm}
\usepackage{myalgorithmic}
\usepackage{algorithm}
\usepackage{url}

\newtheorem{proposition}{Proposition}
\newtheorem{theorem}{Theorem}

\newtheorem{lemma}{Lemma}

\newcommand{\R}{\mathbbm{R}}

\DeclareMathOperator{\Cl}{cl}
\DeclareMathOperator{\Facets}{F}
\DeclareMathOperator{\Vertices}{V}
\DeclareMathOperator{\Sign}{\textsc{sign}}

\newcommand{\card}[1]{\lvert {#1} \rvert}
\newcommand{\order}[1]{\textrm{O}({#1})}%
\newcommand{\setdef}[2]{\left\{{#1}\ :\ {#2}\right\}}
\newcommand{\cl}[1]{\Cl({#1})}
\newcommand{\signvec}[1]{\text{\textsc{#1}}}

\title{Computing the Face Lattice of a Polytope\linebreak from its Vertex-Facet
  Incidences}

\author{Volker Kaibel and Marc E.~Pfetsch%
  \thanks{ TU~Berlin, Fakult\"at~II, Institut f\"ur Mathematik, MA
    6--2,  Stra\ss e des 17. Juni~136, 10623~Berlin,
    \texttt{\{kaibel,pfetsch\}@math.tu-berlin.de}.  The first author has
    been supported by a DFG Gerhard-Hess-Forschungsf\"orderungspreis
    donated to G\"unter M. Ziegler (Zi 475/2-3).}
}

\begin{document}
\maketitle

\begin{abstract}
  We give an algorithm that constructs the Hasse diagram of the face
  lattice of a convex polytope~$P$ from its vertex-facet incidences in
  time~$\order{ \min\{n,m\} \cdot \alpha \cdot \varphi}$, where~$n$ is
  the number of vertices, $m$ is the number of facets, $\alpha$ is the
  number of vertex-facet incidences, and~$\varphi$ is the total number
  of faces of~$P$. This improves results of Fukuda and
  Rosta~\cite{FukR94}, who described an algorithm for enumerating all
  faces of a $d$-polytope in~$\order{\min\{n,m\}\cdot d\cdot
    \varphi^2}$ steps. For simple or simplicial $d$-polytopes our
  algorithm can be specialized to run in time $\order{d \cdot \alpha
    \cdot \varphi}$. Furthermore, applications of the algorithm to
  other atomic lattices are discussed, e.g., to face lattices of
  oriented matroids.
\end{abstract}


\bigskip


\bigskip

\section{Introduction}
\label{sec:intro}

\noindent Let $P$ be a $d$-polytope, i.e., a $d$-dimensional bounded
convex polyhedron. It is well-known that the set $\mathcal{F}$ of its
faces (including~$\varnothing$ and~$P$ itself), ordered by inclusion,
is a graded, atomic, and coatomic lattice: the \emph{face lattice} of
$P$. In particular, each face can be identified with its set of
vertices or the set of facets it is contained in. In this paper, a
face is usually identified with its vertex set. We define $\varphi :=
\card{\mathcal{F}}$ and denote by $\mathcal{L}$ the Hasse diagram (as
an abstract graph) of the face lattice. Hence, $\mathcal{L}$ is a
directed rooted acyclic graph whose nodes correspond to the elements
of $\mathcal{F}$. If $\ell_H, \ell_G$ are nodes in $\mathcal{L}$ and
$H,G \in \mathcal{F}$ are the corresponding faces of~$P$, then there
is an arc~$(\ell_H, \ell_G)$ in $\mathcal{L}$ if and only if $H
\subsetneq G$ and $\dim(G) = \dim(H)+1$.

The \emph{combinatorial face lattice enumeration problem} is the
following: given a vertex-facet incidence matrix of $P$ (see
Section~\ref{sec:algo} for a definition), construct the Hasse diagram
$\mathcal{L}$ of the face lattice. By definition, $\mathcal{L}$ is
unlabeled. Nevertheless, it might be desired to label each node of
$\mathcal{L}$ corresponding to a face~$F$ with the set of (indices of)
vertices contained in~$F$, the set of (indices of) facets
containing~$F$, or with the dimension of~$F$.

Fukuda and Rosta \cite{FukR94} gave an algorithm for the combinatorial
face lattice enumeration problem for $d$-polytopes~$P$ which runs
in~$\order{\min\{n,m\}\cdot d\cdot \varphi^2}$ time, where $m$ is the
number of facets and $n$ is the number of vertices of~$P$. Since
$\varphi$ can be exponential in $n$ and $m$ (consider the $d$-simplex,
for instance) it is desirable to have an algorithm whose running time
depends only linearly on~$\varphi$ (and polynomially on~$n$ and~$m$).
The main result of this paper is such an algorithm.

For the \emph{geometric face lattice enumeration problem}, which asks
for the face lattice of a polytope that is given by an inequality
description, there are algorithms satisfying this requirement on the
running time, e.g., by Fukuda, Liebling, and Margot~\cite{FukLM97}.
However, in our context no geometric data are available.

Ganter~\cite{Gan87} described an algorithm which, given the incidences
of atoms and coatoms of a general atomic and coatomic lattice,
enumerates all elements of the lattice in lexicographic order, where
each element is identified with the set of atoms below it (which are
ordered arbitrarily).  Specialized to our situation, one obtains an
algorithm that computes all vertex sets of faces of~$P$ in
$\order{\min\{n,m\} \cdot \alpha \cdot \varphi}$ steps, where $\alpha$
is the number of vertex-facet incidences of~$P$.  Note that
$d\cdot\max\{n,m\}\leq\alpha\leq n\cdot m$, in particular, $\alpha$ is
bounded polynomially in~$n$ and~$m$. This algorithm, however, does not
compute the inclusion relations between the faces, i.e., the edges of
the Hasse diagram of the face lattice. Of course, once all (vertex
sets of) faces are computed, one may construct the Hasse diagram in an
obvious way after\-wards, but this would require a number of steps which
is quadratic in the total number~$\varphi$ of faces.

Inspired by Ganter's algorithm, we developed the (quite different)
algorithm presented below, which computes the entire Hasse diagram in
the same running time of $\order{ \min\{n,m\} \cdot \alpha \cdot
  \varphi}$, see Theorem~\ref{thm:alg1}. It requires
$\order{\varphi\cdot\min\{n,m\}}$ memory (without output storage). In our
algorithm, the vertex set of each face or the set of facets it is
contained in, as well as its dimension, is readily available (or can
be computed without increasing the asymptotic running time). Of
course, this may increase the (output) storage requirements
significantly.

Fukuda and Rosta~\cite{FukR94} also considered the combinatorial face
lattice enumeration problem for the special case of simple or
simplicial polytopes. They presented an algorithm that computes the
face lattice of a simple polytope in $\order{d\cdot\varphi}$ steps,
provided that in addition to the vertex-facet incidences an acyclic
orientation of the graph of the polytope is given that induces
precisely one sink on every non-empty face. Such an orientation is
called a \emph{good orientation} or an \emph{abstract objective
  function orientation}. Unfortunately, no polynomial time algorithm
is known that computes a good orientation of a simple polytope~$P$ ---
neither if~$P$ is given by its vertex-facet incidences nor if it is
specified by its whole face lattice.

For simple or simplicial polytopes, our algorithm can be specialized
such that it computes the Hasse diagram of the face lattice in
$\order{d\cdot\alpha\cdot\varphi}$ steps from the vertex-facet
incidences, where no good orientation is required (see
Section~\ref{subsec:simp}).

In Section \ref{subsec:skel} we give a rough sketch of the algorithm,
which is followed by a more detailed description in
Sections~\ref{subsec:clos}, \ref{subsec:min}, and \ref{subsec:locate}.
In Section \ref{subsec:analysis} we analyze the algorithm. We present
the specialization of the algorithm for simple or simplicial polytopes
in Section~\ref{subsec:simp} and a variant that computes the
$k$-skeleton in Section~\ref{subsec:kskel}.  Furthermore, in
Section~\ref{subsec:hasse} a version that needs significantly less
memory is described which enumerates just the faces together with
their dimensions (i.e., without the edges of the Hasse diagram).
Finally, a modification that computes the face lattice of an oriented
matroid from its  cocircuits (Section~\ref{subsec:OM}) is
explained.

For the basic properties of polytopes that are important in our
context, we refer to Ziegler's book~\cite{Zie95}. The few concepts
from the theory of algorithms and data structures that play a role in
the paper can be found in any corresponding textbook (e.g. in the one
by Cormen, Leiserson, Rivest, and Stein~\cite{CorLRS01}).
Our running time estimates refer to the uniform time measure (i.e.,
every arithmetic operation/comparison takes one unit of time), while
our statements on memory requirements refer to the bit model.

\section{The Algorithm}
\label{sec:algo}

Define $m$ to be the number of facets and $n$ the number of vertices
of the $d$-polytope~$P$. Let $A = (a_{fv}) \in \{0,1\}^{m \times n}$
be a \emph{vertex-facet incidence matrix} of~$P$. Hence the facets
of~$P$ can be identified with $F := \{1,\dots,m\}$ and its vertices
can be identified with $V := \{1, \dots, n\}$, such that $a_{fv} = 1$
if facet~$f$ contains vertex $v$, and $a_{fv} = 0$ otherwise. Denote
by $\alpha$ the number of vertex-facet incidences, i.e., the number of
ones in $A$. For $S \subseteq V$, define $\Facets(S) := \{ f \in F :
a_{fs} = 1$ for all $s \in S\}$, the set of facets containing all
vertices of $S$. For $T \subseteq F$, define $\Vertices(T) := \{ v \in
V : a_{tv} = 1$ for all $t \in T\}$, the set of vertices contained in
all facets of $T$.

For $S \subseteq V$, the set $\cl{S} := \Vertices(\Facets(S))$
is the (vertex set of) the smallest face of $P$ containing $S$ (in
lattice theoretic terms, the \emph{join} of the elements in~$S$). One
can check easily that this defines a \emph{closure map} on the subsets
of~$V$, i.e., for all $S,S' \subseteq V$ we have:
$$
S \subseteq \cl{S},\qquad
\cl{\cl{S}} = \cl{S},\qquad 
S \subseteq S' \;\; \Rightarrow \;\; \cl{S} \subseteq \cl{S'}.
$$
The faces of $P$ correspond exactly to the \emph{closed sets} of
$V$ with respect to this closure map (i.e., sets $S\subseteq V$ with
$\cl{S}=S$). Our algorithm crucially relies on the fact that closures
can be computed fast (see Section~\ref{subsec:clos}).

\subsection{The Skeleton of the Algorithm}
\label{subsec:skel}

The strategy is to build up the Hasse diagram~$\mathcal{L}$ of the
face lattice from bottom ($\varnothing$) to top ($P$). Consequently,
$\mathcal{L}$ is initialized with the single face~$\varnothing$ and
then enlarged iteratively by adding out-neighbors of nodes that have
already been constructed. We will say that a face has been
\emph{seen}, once its corresponding node in~$\mathcal{L}$ has been
constructed.

During the algorithm, we keep a set $\mathcal{Q}$ containing those
faces that we have seen so far, but for which we have not yet inserted
their out-arcs into the Hasse diagram. At each major step, we remove a
face $H$ from the set~$\mathcal{Q}$ and construct the
set~$\mathcal{G}$ of all faces $G$ with $H\subsetneq G$ and
$\dim(G)=\dim(H)+1$. For each face $G\in\mathcal{G}$ we check whether
it has already been seen. If this is not the case, then a new node
in~$\mathcal{L}$ representing~$G$ is constructed, and~$G$ is added
to~$\mathcal{Q}$. In any case, an arc from the node corresponding
to~$H$ to the node corresponding to~$G$ is inserted
into~$\mathcal{L}$.
 
In order to compute the set~$\mathcal{G}$, we exploit the fact
that~$\mathcal{G}$ consists of the inclusion minimal faces among the
ones that properly contain~$H$.  Since the face lattice of a polytope
is atomic, each face $G\in\mathcal{G}$ must be of the form
$H(v):=\cl{H\cup\{v\}}$ for some vertex (atom)~$v$; in particular, the
Hasse diagram has at most $n \cdot \varphi$ arcs.  Thus, we first
construct the collection~$\mathcal{H}$ of all sets $H(v)$, $v\in
V\setminus H$, and then compute~$\mathcal{G}$ as the set of inclusion
minimal sets of~$\mathcal{H}$.

Computing $H(v)$ for some $v\in V\setminus H$ requires determining a
closure. In Section~\ref{subsec:clos}, we describe a method to perform
this task in $\order{\alpha}$ steps. Determining the inclusion minimal
sets in the collection~$\mathcal{H}$ clearly could be done in
$\order{n^3}$ steps by pairwise comparisons, since $\mathcal{H}$ has
at most $n$ elements, each of size at most~$n$. In
Section~\ref{subsec:min} we show that this can even be performed in
$\order{n^2}$ time.

Another crucial ingredient is a data structure, described in
Section~\ref{subsec:locate}, that allows us to locate the node
in~$\mathcal{L}$ representing a given face~$G$ or to assert that~$G$
has not yet been seen. This can be performed in $\order{\alpha}$
steps.

A summary of the analysis of the time complexity of the algorithm,
along with a pseudo-code description of it, is given in
Section~\ref{subsec:analysis}.

\subsection{Computing Closures}
\label{subsec:clos}

In order to be able to compute closures fast, we store the incidence
matrix~$A$ in a \emph{sorted sparse format} both in a row and column
based way. For each vertex $v \in V$, the elements in
$\Facets(\{v\})\subseteq\{1,\dots,m\}$ are stored increasingly in a
list.  Similarly, for each facet $f \in F$, we store the sorted set
$\Vertices(\{f\})$ in a list. This preprocessing can be performed in
$\order{n \cdot m}$ time (which is dominated by $\order{n \cdot
  \alpha}$ and thus does not influence the estimate of the asymptotic
running time in Proposition~\ref{prop:alg1} below). The sorted sparse
format uses $\order{\alpha\cdot\log\max\{n,m\}}$ storage.

Whenever we want to compute the closure of a set $S \subseteq V$, the
first step is to compute $\Facets(S)$, i.e., the intersection of the
lists $\Facets(\{v\})$, $v\in S$. Since the intersection of two sorted
lists can be computed in time proportional to the sum of the lengths
of the two lists and because the intersection of two lists is at most
as long as the shorter one, $\Facets(S)$ can be computed in time
$\order{\sum_{v \in S} \card{\Facets(\{v\})}}\subseteq\order{\alpha}$.
Similarly, $\Vertices(T)$ can be computed in time $\order{\alpha}$ for
a set $T \subseteq F$.

\begin{lemma}
  \label{lem:clos} 
  The closure $\cl{S}$ of a set $S \subseteq V$ can be computed
  in $\order{\alpha}$ steps (provided that the vertex-facet incidence
  matrix is given in the sorted sparse format).
\end{lemma}

\subsection{Identifying the Minimal Sets}
\label{subsec:min}

Suppose that $H\subsetneq V$ is a face of~$P$ and~$\mathcal{H}$ is the
collection of all faces $H(v)=\cl{H\cup\{v\}}\subseteq V$, $v \in V
\setminus H$.

Our procedure to identify the set~$\mathcal{G}$ of minimal sets in the
collection~$\mathcal{H}$ starts by assigning a label \emph{candidate}
to each vertex in $V\setminus H$. Subsequently, the label
\emph{candidate} of each vertex will either be removed or replaced by
a label \emph{minimal}. We keep the following three invariants: For
each vertex $v$ that is labeled \emph{minimal} we have $H(v) \in
\mathcal{G}$; if two different vertices~$v$ and~$w$ both are labeled
\emph{minimal}, then we have $H(v) \neq H(w)$; $\mathcal{G}$ is
contained in the set of all $H(v)$ for which $v$ is labeled
\emph{minimal} or \emph{candidate}.  Clearly, if no vertex is labeled
\emph{candidate} anymore, the set of vertices labeled \emph{minimal}
is in one-to-one correspondence to~$\mathcal{G}$ via $H(\cdot)$.

Suppose there is still some $v$ labeled \emph{candidate} available. If
$H(v) \setminus \{v\}$ contains some vertex~$w$, then we have $H(w)
\subseteq H(v)$, because $H(w)$ is the intersection of all faces
containing~$H$ and~$w$, and one of these faces is $H(v)$. Hence, if
$w$ is labeled \emph{minimal} or \emph{candidate}, we remove the label
\emph{candidate} from~$v$; otherwise we label~$v$ \emph{minimal}.

It follows by induction that the three invariants are satisfied
throughout the procedure. Moreover, at each major step (choosing a
\emph{candidate}~$v$) the number of \emph{candidate} labels decreases
by one. Since each such step takes $\order{n}$ time, the entire
procedure has complexity $\order{n^2}$.

\begin{lemma}
  \label{lem:min}
  The set~$\mathcal{G}$ of inclusion minimal sets
  in~$\mathcal{H}=\setdef{H(v)}{v\in V\setminus H}$ can be identified
  in $\order{n^2}$ steps.
\end{lemma}

\subsection{Locating Nodes}
\label{subsec:locate}

During the algorithm, we have to keep track of the faces that we have
seen so far and their corresponding nodes in~$\mathcal{L}$.  To this
end, we maintain a special data structure, the \emph{face tree}.  In
this data structure, a face \mbox{$S= \{s_1, \dots, s_k\}\subseteq V$} (with
$s_1 < \dots < s_k$) is represented by the lexicographically smallest
set $C(S) \subseteq S$ that generates~$S$, i.e., $\cl{C(S)} = S$. We
call~$C(S)$ the \emph{canonical spanning set} of the face $S$.  The
map $C(\cdot)$ is one-to-one; its inverse map is the closure map.

The set $C(S)$ can be computed efficiently as follows. For $k=1,2$
set~$C(S) := S$.  For $k \geq 3$, $C(S)$ is computed iteratively:
initialize $C(S)$ with the set $\{s_1,s_2\}$; at each iteration extend
$C(S)$ by the smallest $s_i$ such that $\cl{C(S)} \subsetneq \cl{C(S)
  \cup \{s_i\}}$. Note that $\card{C(S)} \leq \dim(S)+1 \leq d+1$.
Recall that we stored the vertex-facet incidences in the sorted sparse
format (see Section~\ref{subsec:clos}). Similarly to the method for
computing closures, this computation can be performed in
$\order{\alpha}$ steps, since just the intersections $\Facets(\{s_1\})
\cap \dots \cap \Facets(\{s_i\})$, $i=1,\dots,k$, have to be computed
iteratively.  Then, $C(S)$ is obtained as the set of those $s_i$ for
which the intersection becomes smaller.

We now explain the structure of the face tree. Its arcs are directed
away from the root. They are labeled with vertex numbers, such that no
two arcs leaving the same node have the same label and on every
directed path in the tree the labels are increasing. Via the sets of
labels on the paths from the root, the nodes of the tree correspond to
the sorted sets $C(S)$ for the faces~$S \subseteq V$ that have been
seen so far. In particular, the root node represents the face
$\varnothing$. Each node has a pointer to the corresponding node of
$\mathcal{L}$.  By construction, the depth of the tree is bounded by
$d+1$.

Suppose we want to find the node $\ell_S$ corresponding to some face
$S\subseteq V$ in the part of~$\mathcal{L}$ that we have already
constructed or to assert that this face has not yet been seen. We
first sort~$S$ (a subset of $\{1,\dots,n\}$) increasingly
in~$\order{n}$ steps (by counting or bucket sort,
see~\cite[Chap.~8]{CorLRS01}) and compute~$C(S)$ in $\order{\alpha}$
steps. Then, starting from the root, we proceed (as long as possible)
downwards in the face tree along arcs labeled by the successive
elements of $C(S)$.  Either we find an existing node in the tree which
corresponds to~$S$, or we have to introduce new labeled arcs (and
nodes) into the tree until we have constructed a node
representing~$S$.

In the latter case, it might be necessary to construct an entire new
path in the tree. The definition of the canonical spanning sets $C(S)$
ensures that all ``intermediate nodes'' on that path will correspond
to canonical spanning sets of faces as well. Hence, the number of
nodes in the face tree always will be bounded by~$\varphi$, the total
number of faces of the polytope. The faces represented by intermediate
nodes will be seen later in the algorithm. Consequently, the
corresponding pointers to~$\mathcal{L}$ are set to \texttt{nil} for
the meantime. Later in the algorithm, when we are searching for the
face represented by such a tree-node for the first time, the
\texttt{nil}-pointer will indicate that this face is not yet
represented in~$\mathcal{L}$. The \texttt{nil}-pointer is then
replaced by a pointer to a newly created node representing the face
in~$\mathcal{L}$.

In any case, since the face tree has depth at most $d+1$ and the
out-degree of each node is at most $n$, we need a total time of
$\order{n + \alpha + (d+1) \cdot n} = \order{\alpha}$ to either locate
or create the tree-node representing a certain face.

\begin{lemma}
  \label{lem:locate}
  Using the face tree, it is possible to locate or create the node
  in~$\mathcal{L}$ representing a face in $\order{\alpha}$ steps
  (provided the vertex-facet incidence matrix is stored in the sorted
  sparse format).
\end{lemma}

In the description given above, we have assumed that for each node in
the face tree the out-arcs are stored in a list which is searched
linearly for a certain label when walking down the tree. Of course,
one can store the set of out-arcs in a balanced search tree (see,
e.g., \cite[Chap.~13]{CorLRS01}), allowing to perform the search for a
certain label in logarithmic time.  After computing~$C(S)$ for a
face~$S$ (in $\order{\alpha}$ time), this allows to locate or create
the node corresponding to~$S$ in the face tree in $\order{(d+1) \cdot
  \log n}$ steps. The total running time remains $\order{\alpha}$;
nevertheless this might speed up the algorithm in practice.

Instead of using the face tree, one can also store the faces in a
balanced search tree. Again, the faces are represented by their
canonical spanning sets, which are ordered lexicographically.  Once
$C(S)$ is computed for a face~$S$, searching~$S$ can be performed in
$\order{(d+1)\cdot\log\varphi}\subseteq\order{(d+1)\cdot\min\{n,m\}}$
steps (since $\varphi\leq 2^{\min\{n,m\}}$). This yields the same
total asymptotic running time, but searching the tree takes more (or
the same) time compared to the variant of the face tree with balanced
search trees at its nodes, since $\log n\le\min\{n,m\}$.

\subsection{The Analysis}
\label{subsec:analysis}

We summarize the algorithm in pseudo-code
(Algorithm~\ref{algo:enumeration}):

\begin{algorithm}
\caption{Computing the face lattice of a polytope from its incidences}
\label{algo:enumeration}
\begin{algorithmic}[1]
  \STATE \textbf{Input:} incidence matrix of a polytope~$P$
  \STATE \textbf{Output:} Hasse diagram~$\mathcal{L}$ of the face lattice of~$P$
  \STATE\label{line:init1} initialize $\mathcal{L}$ and a face tree with
  $\ell_{\varnothing}$ corresponding to the empty face
  \STATE\label{line:init2} initialize a set~$\mathcal{Q} \subseteq
  \{$nodes of $\mathcal{L}\} \times \{$subsets of $V\}$ by $(\ell_{\varnothing},\varnothing)$
  \WHILE{$\mathcal{Q} \neq \varnothing$}
    \STATE\label{line:choose} choose some $(\ell_H,H) \in \mathcal{Q}$ and remove it from~$\mathcal{Q}$
    \STATE\label{line:H} compute the collection~$\mathcal{H}$ of all $H(v)$, $v \in V\setminus H$
    \STATE\label{line:G} compute the set~$\mathcal{G}$ of minimal sets 
    in~$\mathcal{H}$
    \FOR{each $G\in\mathcal{G}$}
      \STATE\label{line:locateG} locate/create the node~$\ell_G$ corresponding
      to~$G$ in~$\mathcal{L}$ 
      \IF{$\ell_G$ was newly created}
        \STATE add~$(\ell_G,G)$ to~$\mathcal{Q}$
      \ENDIF
      \STATE add the arc $(\ell_H,\ell_G)$ to~$\mathcal{L}$
    \ENDFOR
  \ENDWHILE
\end{algorithmic}
\end{algorithm}

\medskip

\begin{proposition} \label{prop:alg1}
  Algorithm~\ref{algo:enumeration} computes the Hasse diagram of the
  face lattice of a polytope~$P$ from its vertex-facet incidences in
  $\order{n \cdot \alpha \cdot \varphi}$ time. It can be implemented
  such that its space requirements (without output space) are bounded
  by~$\order{\varphi\cdot n}$.
\end{proposition}

\begin{proof}
  Algorithm \ref{algo:enumeration} works correctly by the discussion
  above.
  
  Step~\ref{line:H} can be performed in $\order{n\cdot\alpha}$ steps
  by Lemma~\ref{lem:clos}. Lemma~\ref{lem:min} shows that we can
  execute Step~\ref{line:G} in~$\order{n^2} \subseteq \order{n \cdot
    \alpha}$ time. Hence, Steps~\ref{line:H} and~\ref{line:G} in total
  contribute at most $\order{n \cdot \alpha \cdot \varphi}$ to the
  running time (since the while-loop is executed once per face).
  
  The body of the for-loop has to be executed for each of the
  $\order{n \cdot \varphi}$ arcs in the Hasse diagram $\mathcal{L}$.
  Lemma~\ref{lem:locate} implies that each execution of the body of
  the for-loop can be performed in $\order{\alpha}$ steps. Thus, the
  claim on the running time follows.
  
  Since each node of the face tree corresponds to a face of $P$, the
  face tree has $\order{\varphi}$ nodes. Each label on an edge of the
  face tree needs at most~$\order{\log n}$ bits, and we can estimate
  the space requirements of any of the (internal and external)
  pointers by $\order{\log\varphi}\subseteq\order{\min\{n,m\}}$.
  Hence, the face tree needs no more than
  $\order{\varphi\cdot\min\{n,m\}}$ bits.
  
  The space required for the storage of~$\mathcal{Q}$ is bounded by
  $\order{\varphi\cdot n}$, if for each pair
  $(\ell_H,H)\in\mathcal{Q}$ the set~$H$ is stored as a \emph{bit
    set}, i.e., the characteristic vector of $H\subseteq V$ is stored
  bit by bit.
\end{proof}

If $m<n$, then it is more efficient to apply Algorithm
\ref{algo:enumeration} to the incidences of the dual polytope, i.e.,
to the transpose of the incidence matrix. Of course, after the
computations the roles of vertices and facets have to be exchanged
again. This yields the main result of the paper.

\begin{theorem} \label{thm:alg1}
  The Hasse diagram of the face lattice of a polytope~$P$ can be
  computed from the vertex-facet incidences of~$P$ in
  $\order{\min\{n,m\} \cdot \alpha \cdot \varphi}$ time, where~$n$ is
  the number of vertices, $m$ is the number of facets, $\alpha$ is the
  number of vertex-facet incidences, and~$\varphi$ is the total number
  of faces of~$P$. The space requirements of the algorithm (without
  output space) can be bounded by~$\order{\varphi\cdot\min\{n,m\}}$.
\end{theorem}

Whenever a new node representing a face~$G$ in the Hasse
diagram~$\mathcal{L}$ is constructed, we can label that node with the
vertex set of~$G$, the set of facets containing~$G$, or with the
dimension of~$G$ without (asymptotically) increasing the running time
of the algorithm. The output, however, might become much larger with
such labelings. For instance, labeling the Hasse diagram of the
$d$-cube by vertex labels requires $\Omega(4^d\cdot d)$ output storage
space, while the Hasse diagram with facet labels needs only
$\order{d \cdot 3^d\cdot\log d}$ space.

In practice, the computation can be speeded up by exploiting that
every vertex that is contained in a face~$G$ with $H\subsetneq G$ and
$\dim{G}=\dim{H}+1$ must be contained in some facet which
contains~$H$. Thus, it suffices to consider only the sets $H(v)$, $v
\in \bigl( \bigcup_{f \in \Facets(H)} \Vertices(\{f\}) \bigr)
\setminus H$ in Step~\ref{line:H}.

\section{Extensions}
\label{sec:ext}

\subsection{Simple or Simplicial Polytopes}
\label{subsec:simp}

For a simple $d$-polytope~$P$ with~$n$ vertices, the above procedure
can be implemented to run more efficiently. We have $\alpha = n \cdot d$ in
this case. From the incidences (stored in the sorted sparse format),
the graph $G(P)$ of~$P$ (i.e., all one-dimensional faces) can be
computed in time $\order{n^2 \cdot d}$, since a pair of vertices forms
an edge if and only if it is contained in $d-1$ common facets.

After initialization with the vertices instead of $\varnothing$ (in
Steps \ref{line:init1} and \ref{line:init2}), Steps~\ref{line:H}
and~\ref{line:G} can now be simplified. Consider an arbitrary
vertex~$w \in H$. For each neighbor $v \notin H$ of $w$ in $G(P)$,
$H(v)$ yields the other end node of an arc in the Hasse diagram; and
each out-arc of $H$ is produced this way.  Thus, we can avoid
constructing non-minimal faces in Step~\ref{line:H}. Hence,
Step~\ref{line:G} can be skipped.  The total running time for simple
$d$-polytopes decreases to~$\order{d \cdot \alpha \cdot \varphi}$
(since the body of the for-loop is executed at most $d\cdot\varphi$
times). 

The space complexity stays $\order{\varphi\cdot n}$ (see
Proposition~\ref{prop:alg1}). It can, however, be reduced to
$\order{\varphi\cdot m}$ (we have $m\le n$ for simple polytopes):
instead of storing pairs $(\ell_H,H)$ in the set~$\mathcal{Q}$, we
store the pairs $(\ell_H,\Facets(H))$, since~$\card{\Facets(H)}\le m$.
Converting between $H$ and $\Facets(H)$ can be performed
in~$\order{\alpha}$ steps and hence does not increase the asymptotic
total running time.

By duality, the same running times and space requirements can be
achieved for simplicial polytopes.

Similarly to the situation with general polytopes, for both simple and
simplicial polytopes we can also output for each face its vertices,
the facets containing it, or its dimension without (asymptotically)
increasing the running time.

\subsection{The {\it \textbf{k}}-Skeleton}
\label{subsec:kskel}

A variant of Algorithm~\ref{algo:enumeration} computes the Hasse
diagram of the $k$-skeleton (all faces of dimension at most $k$) of a
polytope $P$. One simply prevents the computation of faces of
dimensions larger than~$k$ by not inserting any $(k-1)$-face into the
list~$\mathcal{Q}$.  This leads to an $\order{n \cdot \alpha \cdot
  \varphi^{\leq k}}$ time algorithm, where $\varphi^{\leq k}$ is the
number of faces of~$P$ of dimension at most $k$.

\subsection{Computing the ``Hasse Diagram without Edges''}
\label{subsec:hasse}

If we only want to compute the faces of $P$ together with their
descriptions and dimensions (i.e., the ``Hasse diagram without
edges''), there exists a variant of Algorithm~\ref{algo:enumeration}
with the same asymptotic running time, but sig\-nifi\-cantly reduced
space requirements.  The difference is that no face tree is used, and
the set $\mathcal{Q}$ is organized as a stack, i.e., the faces are
investigated in a depth-first search manner. At each step, we take a
face $H$ from the stack, output it, and compute the set $\mathcal{G}$
of $(\dim H + 1)$-faces containing $H$, like in Steps \ref{line:H} and
\ref{line:G} of Algorithm~\ref{algo:enumeration}.  This needs time
$\order{n \cdot \alpha}$ for each $H$. The for-loop beginning at Step
\ref{line:locateG}, including the search in the face tree, is replaced
by an efficient way to decide which of the faces in $\mathcal{G}$ is
put onto the stack $\mathcal{Q}$, such that every face appears exactly
once on the stack during the algorithm. For this, we compute for each
face $G \in \mathcal{G}$ a unique canonical facet $H'$ of it. We
put~$G$ onto the stack if and only if $H = H'$. This ensures that each
face is picked exactly once.

We take $H'$ as the closure of a set $D(G)$, which is computed similar
to the set $C(G)$ of Section \ref{subsec:locate}, except that we
reject vertices which would produce $G$. More precisely, let $G =
\{g_1, g_2, \dots, g_l\}$, with $g_1 < g_2 < \dots <g_l$. Initialize
$D(G)$ with $\varnothing$ and in each iteration extend $D(G)$ by the
smallest~$g_i$ such that $\cl{D(G)} \subsetneq \cl{D(G) \cup \{g_i\}}$
and $\cl{D(G) \cup \{g_i\}} \neq G$. After the computation, $H'$, the
closure of $D(G)$, clearly is a proper face of $G$. Moreover, it is a
facet of~$G$, since otherwise there exists a vertex $g \in G\setminus
H'$, such that $\cl{H' \cup \{g\}} \subsetneq G$. But then $g$ would
have been included into~$D(G)$ when it was considered.  Hence, $D(G)$
is the lexicographically smallest subset of $G$ which spans a facet
of~$G$.  It can be computed in time $\order{\alpha}$, and hence,
checking for all faces $G\in\mathcal{G}$ whether~$H$ is the canonical
facet $D(G)$ of $G$ can be performed in $\order{n \cdot \alpha}$ time.

Altogether, this leads to an $\order{n \cdot \alpha \cdot \varphi}$
time algorithm (see the proof of Proposition~\ref{prop:alg1}). The
algorithm needs $\order{n^2\cdot d\cdot\log n }$ space
for~$\mathcal{Q}$; since the depth of $\mathcal{Q}$ is at most $d+1$,
there are never more than $n\cdot (d+1)$ sets on the stack, each of
size at most~$n$.  Additionally, we need
$\order{\alpha\cdot\log\max\{n,m\}}$ space for storing the incidences
in the sorted sparse format.  Applying this method to the dual
polytope, if necessary, we obtain an $\order{\min\{n,m\} \cdot \alpha
  \cdot \varphi}$ time algorithm.

\subsection{Oriented Matroids}
\label{subsec:OM}

Algorithm~\ref{algo:enumeration} can be used for the enumeration of
the elements of any atomic lattice provided a subroutine is available
that computes the join of a set of atoms. For instance, this holds for
every atomic and coatomic lattice if the atom-coatom incidences are
given, because in this case one can compute the joins of atoms
similarly to the case of face lattices of polytopes.

In the following, we describe such an application of our algorithm to
oriented matroids. The set of covectors of an oriented matroid with
ground set $\{1,\dots,k\}$ is a subset of $\{-,0,+\}^k$ that satisfies
certain axioms. We refer to Bj\"orner, Las~Vergnas, Sturmfels, White,
and Ziegler \cite[Chap.~4]{BjoLSWZ99} for the definitions and concepts
that are relevant in the following.  A specific, but illustrative,
example arises from any finite set $X$ of points in~$\R^d$ as follows.
For every linear functional $\varphi\in(\R^d)^{\star}$ denote by
$\Sign(\varphi)\in\{-,0,+\}^X$ the vector whose component
corresponding to $x\in X$ encodes the sign of~$\varphi(x)$.
Then~$\setdef{\Sign(\varphi)}{\varphi\in(\R^d)^{\star}}$ is the set of
covectors of an oriented mat\-roid~$\mathcal{O}(X)$.

For $\signvec{v},\signvec{w}\in\{-,0,+\}^k$ the \emph{separation set}
of~$\signvec{v}$ and~$\signvec{w}$
contains all indices~$i$ such that one of $\signvec{v}_i$,
$\signvec{w}_i$ is~$+$, and the other one is~$-$. The
\emph{composition}~$\signvec{v}\circ\signvec{w}$ of~$\signvec{v}$
and~$\signvec{w}$ is defined by
$(\signvec{v}\circ\signvec{w})_i:=\signvec{v}_i$ if $\signvec{v}_i\neq
0$ and $(\signvec{v}\circ\signvec{w})_i:=\signvec{w}_i$ otherwise.

We define a partial order~$\preceq$ on $\{-,0,+\}^k$, where
$\signvec{v}\preceq\signvec{w}$ holds if and only if for all~$i$ we
have $\signvec{v}_i=0$ or $\signvec{v}_i=\signvec{w}_i$. The
$\preceq$-minimal elements among the nonzero covectors of an oriented
matroid are called its \emph{cocircuits}. If one adjoins an
artificial maximal element~$\hat{1}$ to the poset formed by the
covectors of an oriented matroid (ordered by~$\preceq$),
then one obtains its (\emph{big}) \emph{face lattice}.

If, in the above example, $X$ is the vertex set of a
polytope~$P\subset\R^d$, then the faces of~$P$ correspond to the
\emph{positive covectors} (i.e., the covectors with no component equal
to $-$) of $\mathcal{O}(X)$. The facets of~$P$ correspond to the
positive cocircuits of $\mathcal{O}(X)$.  The face lattice of~$P$ is
anti-isomorphic to a sublattice of the face lattice of
$\mathcal{O}(X)$.

The face lattice of an oriented matroid is atomic and coatomic; its
atoms are the cocircuits, and its coatoms are called \emph{topes}.
Hence, we can compute its Hasse diagram from the abstract atom-coatom
incidences as above.

However, this is not the usual way to encode an oriented matroid. It
is more common to specify an oriented matroid by its cocircuits.  The
join of two covectors simply is their composition, if their separation
set is empty, or $\hat{1}$ otherwise.  Such a composition can be
computed in $\order{k}$ steps, which enables us to compute the face
lattice (efficiently) from its cocircuits by a variant of
Algorithm~\ref{algo:enumeration}. 

In Step~\ref{line:choose}, $H$ now is a face of the oriented matroid,
i.e., a covector. In Step~\ref{line:H}, one has to compute the joins
of~$H$ with every cocircuit~$\signvec{v}\not\preceq H$. Thus,
Step~\ref{line:H} can be performed in $\order{n \cdot k \cdot
  \varphi}$ steps altogether (where~$\varphi$ is the total number of
covec\-tors and~$n$ is the number of cocircuits).  We
do not know any  method to perform Step~\ref{line:G} faster than
by pairwise comparisons, which take $\order{n^2\cdot k \cdot \varphi}$
time in total.

The face tree is organized similarly to the description in
Section~\ref{subsec:locate}.  One fixes an
ordering~$\signvec{C}_1,\dots,\signvec{C}_n$ of the cocircuits. For a
covector~$\signvec{S}$ let $\{i_1,\dots,i_r\}$ ($i_1<\dots<i_r$) be
the index set of cocircuits $\signvec{C}_{i_j}\preceq \signvec{S}$.  Then we
iteratively form the joins of
$\signvec{C}_{i_1},\dots,\signvec{C}_{i_r}$, and let $C(\signvec{S})$
consist of all those indices for which the ``joins change.''
Computing $C(\signvec{S})$ from~$\signvec{S}$ takes $\order{n\cdot k}$
steps. Note that~$\card{C(\signvec{S})}\leq k$.

Using this modified face tree, a given covector~$\signvec{S}$ can now
be searched in the same way as in the case of face lattices of
polytopes.  The depth of the face tree is bounded by~$k$.  Hence,
location/creation of a covector can be performed in $\order{n\cdot k}$
time. The rest of the analysis is similar to the proof of
Proposition~\ref{prop:alg1}.  Thus, by this variant of
Algorithm~\ref{algo:enumeration}, the Hasse diagram of the face
lattice of an oriented matroid can be computed in $\order{n^2\cdot
  k\cdot \varphi}$ steps, requiring $\order{\varphi\cdot k}$ working
space (since $\varphi\le 3^k$).

Finschi~\cite{Fin01} describes a different algorithm that computes the
covectors of an oriented matroid from its cocircuits in $\order{n\cdot
  k^2\cdot \varphi}$ time. His algorithm, however, does not produce
the edges of the Hasse diagram.

The case where the topes (i.e., the $\preceq$-maximal covectors) of an
oriented matroid are given is a bit different. Here, the number of
faces is bounded by $m^2$, where $m$ is the number of topes.  Hence,
the size of the face lattice is polynomial in $m$.  Fukuda, Saito, and
Tamura \cite{FukST91} give an $\order{k^3 \cdot m^2}$ time algorithm
for constructing the face lattice from the maximal covectors.

\section*{Acknowledgements}
We are indepted to Michael Joswig and J\"org Rambau for stimulating
discussions as well as to G\"unter M.~Ziegler for valuable comments on
the paper. We also thank the two referees for their helpful comments.

\end{document}